\documentclass[12pt,a4paper]{article}
\usepackage{latexsym}
\usepackage{amssymb}
\newtheorem{lem}{Lemma}[section]
\newtheorem{prop}{Proposition}[section]
\newtheorem{thm}{Theorem}[section]
\newenvironment{pf}{\textbf{Proof\ }}{\hfill$\Box$\smallskip}
\newcommand{\sub}{\subseteq}
\newcommand{\id}[1]{\left<#1\right>}
\newcommand{\leg}[2]{\left(\frac{#1}{#2}\right)}
\newcommand{\op}{\oplus}
\newcommand{\al}{\alpha}
\newcommand{\ga}{\gamma}
\newcommand{\de}{\delta}
\newcommand{\ze}{\zeta}
\newcommand{\La}{\Lambda}
\newcommand{\Om}{\Omega}

\newcommand{\be}{\mathbf{e}}
\newcommand{\bu}{\mathbf{u}}
\newcommand{\bv}{\mathbf{v}}
\newcommand{\bw}{\mathbf{w}}
\newcommand{\ca}{\mathcal{A}}
\newcommand{\cc}{\mathcal{C}}
\newcommand{\cd}{\mathcal{D}}
\newcommand{\ci}{\mathcal{I}}
\newcommand{\cj}{\mathcal{J}}
\newcommand{\co}{\mathcal{O}}
\newcommand{\cp}{\mathcal{P}}
\newcommand{\cq}{\mathcal{Q}}
\newcommand{\cs}{\mathcal{S}}
\newcommand{\cw}{\mathcal{W}}
\newcommand{\C}{\mathbf{C}}
\newcommand{\F}{\mathbf{F}}
\newcommand{\Q}{\mathbf{Q}}
\newcommand{\R}{\mathbf{R}}
\newcommand{\Z}{\mathbf{Z}}
\title{Conference Matrices and Unimodular Lattices}
\author{Robin Chapman\\
School of Mathematical Sciences\\ University of Exeter\\
Exeter, EX4 4QE, UK\\ \texttt{rjc@maths.ex.ac.uk}}
\date{19 July 2000}
\begin{document}
\maketitle

\section{Introduction}

We use conference matrices to define an action of the complex numbers
on the real Euclidean vector space $\R^n$. In certain cases, the lattice
$D_n^+$ becomes a module over a ring of quadratic integers. We can then
obtain new unimodular lattices, essentially by multiplying the
lattice $D_n^+$ by a non-principal ideal in this ring. We show that
lattices constructed via quadratic residue codes, including the Leech lattice,
can be constructed in this way.

Recall that a \emph{lattice} $\La$ is a discrete subgroup of a finite
dimensional real vector space~$V$. We suppose that $V$ has a given
Euclidean inner product $(\bu,\bv)\mapsto\bu\cdot\bv$
and the rank of $\La$ equals the dimension of~$V$.
In this case $\La$ has a bounded
fundamental region in~$V$. We call the volume of such a fundamental
region (measured with respect to the Euclidean structure on $V$)
the \emph{volume} of the lattice~$\La$.

The lattice $\La$ is \emph{integral} if $\bu\cdot\bv\in\Z$
for all $\bu$, $\bv\in\La$. It is \emph{even} if $|\bu|^2=\bu\cdot\bu\in2\Z$
for all $\bu\in\La$. Even lattices are necessarily integral. The
lattice $\La$ is \emph{unimodular} if $\La$ is integral and has volume~1.
It is well known \cite[Chapter~VIII, Theorem~8]{Serre}
that if $\La$ is an even unimodular
then the rank of $\La$ is divisible by~8.

For convenience we call the square of the length of a vector its
\emph{norm}. The \emph{minimum norm} of a lattice is the
smallest non-zero norm of its vectors.

\section{Conference matrices}

Let $l$ be a positive integer. A \emph{conference matrix} of order $n$
\cite[Chapter 18]{vLW}
is an $n$-by-$n$ matrix $W$ satisfying
\begin{description}
\item{(a)} the diagonal entries of $W$ vanish, while its off-diagonal entries
lie in $\{-1,1\}$,
\item{(b)} $WW^\top=(n-1)I$, where $I$ denotes the $n$-by-$n$
identity matrix.
\end{description}
Let $\cw_n$ denote the set of skew-symmetric
conference matrices of order~$n$.

Let $W\in\cw_n$. Then $H=I+W$ satisfies
$HH^t=(I+W)(I-W)=I-W^2=I+WW^\top=nI$. As all the entries of $H$ lie
in $\{-1,1\}$ then $H$ is a Hadamard matrix. Consequently
\cite[Theorem 18.1]{vLW} $n=1$, 2 or is a multiple of~4.

Suppose that $n$ is a multiple of 4 and let $l=n-1$.
Fix $W\in\cw_n$ and let  $V=\R^n$ denote the $n$-dimensional real vector
space under the standard Euclidean dot product. Then, since $W^2=-lI$, $V$
becomes also a complex vector space when we define
$$(r+s\sqrt{-l})\bv=\bv(r+sW)$$
for $r$, $s\in\R$. Let $|\bv|=\sqrt{\bv\cdot\bv}$ denote the
Euclidean length of a vector~$\bv\in\R^n$. This action of $\C$ on
$\R^n$ transforms lengths in the obvious way. Let $z^*$ denote
the complex conjugate of the complex number~$z$.

\begin{lem}\label{length}
\begin{description}
\item{(a)}
If $z_1$, $z_2\in\C$ and $\bv_1$, $\bv_2\in\R^n$ then
$(z_1\bv_1)\cdot(z_2\bv_2)=(z_1z_2^*\bv_1)\cdot\bv_2$
\item{(b)}
If $z\in\C$ and $\bv\in\R^n$ then $|z\bv|=|z||\bv|$.
\end{description}
\end{lem}
\begin{pf}
Let $z_j=r_j+s_j\sqrt{-l}$ with $r_j$, $s_j\in\R$. Then
\begin{eqnarray*}
(z_1\bv_1)\cdot(z_2\bv_2)&=&(z_1\bv_1)(z_2\bv_2)^\top\\
&=&\bv_1(r_1I+s_1W)(r_2I+s_2W)^\top\bv_2^\top\\
&=&\bv(r_1I+s_1W)(r_2I-s_2W)\bv^\top\\
&=&\bv((r_1r_2+ls_1s_2)I+(s_1r_2-r_1s_2)W)\bv^\top\\
&=&(z_1z_2^*\bv_1)\cdot\bv_2
\end{eqnarray*}
as claimed.

Consequently
$$|z\bv|^2=(z\bv)\cdot(z\bv)=(zz^*\bv)\cdot\bv=|z|^2\bv\cdot\bv
=|z|^2|\bv|^2.$$
\end{pf}

Thus for fixed nonzero~$z$, the map $\bv\mapsto z\bv$
is a similarity of $\R^n$ with scale factor~$|z|$.

\section{Quadratic fields}

We retain the previous notation.
Suppose in addition that $l=n-1$ is squarefree. Let $K$ denote the
quadratic field $\Q(\sqrt{-l})$. Since $l$ is square-free,
the ring of integers of $K$ is
$$\co=\Z\left[\frac{1+\sqrt{-l}}{2}\right]
=\left\{\frac{a+b\sqrt{-l}}{2}:a,b\in\Z,a\equiv b\pmod{2}\right\}.$$
In particular $\co$ is a Dedekind domain.
We shall show that some of the familiar lattices
in $\R^{l+1}$ are modules for the ring~$\co$.

Let
$$L_0=\{(a_1,\ldots,a_n)\in\Z^n:a_1+\cdots+a_n\equiv 0\pmod{2}\}$$
be the $D_n$ root lattice.

\begin{lem}\label{lemmL0}
The lattice $L_0$ is an $\co$-module.
\end{lem}
\begin{pf}
It suffices to show that $\frac12(1+\sqrt{-l})\bv=\frac12\bv(I+W)\in L_0$
whenever $\bv\in L_0$. Indeed it suffices to show this whenever
$\bv$ lies in a generating set for~$L_0$. Now $L_0$ is generated by the
vectors $2\be_j$ and $\be_j+\be_k$ (for $j\ne k$)
where $\be_j$ denotes the $j$-th unit vector. Firstly $\be_j(I+W)$
is a row of the Hadamard matrix $I+W$. As it contains $n$ instances of
$\pm1$ and $n$ is even, it lies in~$L_0$. Next $\frac12(\be_j+\be_k)(I+W)$
is the sum of two rows of the Hadamard matrix $I+W$. Two rows
of an $n$-by-$n$ Hadamard matrix agree in exactly $n/2$ places.
Hence $\frac12(\be_j+\be_k)(I+W)$ has $n/2$ zeros and
$n/2$ instances of $\pm1$. As $n/2$ is even then
$\frac12(\be_j+\be_k)(I+W)\in L_0$. This completes the proof.
\end{pf}

Now consider the set
$$\cs=\{(a_1,\ldots,a_n):a_j\in\{-1/2,1/2\}\}.$$
The difference of two vectors in $\cs$ lies in $L_0$ if and only
if those vectors agree in an even number of places. Thus there
are exactly two cosets $\bv+L_0$ as $\bv$ runs through~$\cs$.

For each $j$, $\frac12\be_j(I+W)\in\cs$, and for each $j$ and $k$,
$\frac12(\be_j-\be_k)(I+W)\in L_0$ by Lemma~\ref{lemmL0}. Thus the
cosets $\frac12\be_j(I+W)+L_0$ are identical. Let
$$\cs_+=\{\bv\in\cs:\bv-\textstyle{\frac12}\be_1(I+W)\in L_0\}$$
and
$$\cs_-=\cs\setminus\cs_+.$$
As $\frac12\be_j(I+W)-\frac12\be_j(-I+W)=\be_j\notin L_0$ then
$\frac12\be_j(-I+W)\in\cs_-$ for each~$j$.

If $\bv\in\cs$ then $2\bv$ has $n$ entries $\pm1$ and so $2\bv\in L_0$.
It follows that $L_0\cup(\bv+L_0)$ is a lattice, which depends
only on whether $\bv\in\cs_+$ or $\bv\in\cs_-$. We write
$L_+$ for $L_0\cup(\bv+L_0)$ when $\bv\in\cs_+$ and
$L_-$ for $L_0\cup(\bv+L_0)$ when $\bv\in\cs_-$. Both
$L_+$ and $L_-$ are isometric to the lattice usually denoted
by $D_n^+$ \cite[Chapter~4, \S7.3]{CS}. The lattice $D_n^+$ is unimodular
for each $n$ divisible by~4, and it is even unimodular whenever
$n$ is divisible by~8.

\begin{lem}\label{lemmLpm}
If $n$ is divisible by 8 then the lattices $L_+$ and $L_-$ are $\co$-modules.
\end{lem}
\begin{pf}
Let $L=L_+$ or $L_-$. Then $L=L_0+(\bv+L_0)$ for some $\bv\in\cs$
and by Lemma~\ref{lemmL0}
it suffices to show that $\frac12(1+\sqrt{-l})\bv=\frac12\bv(I+W)\in L$.
Note that $\frac14(l+1)$ is an even integer by the hypothesis.

We may assume that $\bv=\frac12\be_1(\pm I+W)$. If $\bv=\frac12\be_1(I+W)$
then
$$\frac12\bv(I+W)=\frac14\be_1(I+W)^2=\frac14\be_1((1-l)I+2W)
=\frac12\be_1(I+W)-\frac{l+1}{4}\be_1$$
which lies in $L$ as $\frac12\be_1(I+W)\in L$.

If
$\bv=\frac12\be_1(-I+W)$
then
$$\frac12\bv(I+W)=\frac14\be_1(-I+W)(I+W)=\frac14\be_1(-(l+1)I)$$
which lies in~$L$.
\end{pf}

Let $\ci$ be an ideal of $\co$. If $M$ is a $\co$-module, then
$\ci M$, defined as the subgroup of $M$ generated by the $\al m$
for $\al\in\ci$ and $m\in M$, is also a $\co$-module.

\begin{thm}\label{main}
Suppose that $l\equiv7$ (mod 8) and that $\ci$ is a nonzero
ideal of $\co$ with norm $N=N(\ci)$. If $L=L_+$ or $L_-$ then
$$L[\ci]=\frac{1}{\sqrt{N}}\ci L$$
is an even unimodular lattice.
Also if $\ci$ and $\cj$ lie in the same ideal class of~$\co$, the lattices
$L[\ci]$ and $L[\cj]$ are isometric.
\end{thm}
\begin{pf}
First of all we show that the index $|L:\ci L|$ equals~$N^{n/2}$.
As an $\co$-module, $L$ is finitely generated. Also if $\al\in\co$
and $\bv\in L$ are nonzero, then $|\al\bv|=|\al||\bv|\ne0$ by
Lemma~\ref{length} and so $L$ is torsion free as an $\co$-module.

By the theory of modules over Dedekind domains \cite[\S 9.6]{PCohn},
as $L$ is a finitely generated torsion-free module over the
Dedekind domain~$\co$, then
$L=L_1\op\cdots\op L_k$ where each $L_j$ is isomorphic to a nonzero ideal
$\ca_j$ of~$\co$. Each of the $\ca_j$ is a free abelian group of rank~2,
and as $L$ is a free abelian group of rank $n$ it follows that $k=n/2$.
Then $\ci L=\ci L_1\op\cdots\op \ci L_{n/2}$
and so $|L:\ci L|=\prod_{j=1}^{n/2}|L_j:\ci L_j|$. But
$$|L_j:\ci L_j|=|\ca_j:\ci\ca_j|=\frac{|\co:\ci\ca_j|}{|\co:\ca_j|}
=\frac{N(\ci\ca_j)}{N(\ca_j)}.$$
But $N(\ci\ca_j)=N(\ci)N(\ca_j)$ and so $|L_j:\ci L_j|=N(\ci)=N$.
Consequently $|L:\ci L|=N^{n/2}$ as claimed.

We now show that $L[\ci]$ is unimodular. The lattice
$\ci L$ is generated by elements $\bu=\al\bv$ where $\al\in\ci$ and
$\bv\in L$. Let $\bu_j=\al_j\bv_j$ ($j=1$, 2) with $\al_j\in\ci$
and $\bv_j\in L$. Then by Lemma~\ref{length},
$$\bu_1\cdot\bu_2=(\al_1\bv_1)\cdot(\al_2\bv_2)
=(\al_1\al_2^*\bv_1)\cdot\bv_2.$$
But $\al_1\al_2^*\in\ci\ci^*=N(\ci)\co$ \cite[\S VIII.1]{HCohn} so that.
$$\bu_1\cdot\bu_2=N(\ga\bv_1)\cdot\bv_2$$
where $\ga\in\co$. As $\ga\bv_1\in L$ (by Lemma~\ref{lemmLpm})
and $L$ is an integral lattice, then $\bu_i\cdot\bu_2\equiv0$
(mod~$N$). Consequently $L[\ci]=N^{-1/2}\ci L$ is an integral
lattice. But $L$ is unimodular, so it has volume~1.
Thus $\ci L$ has volume $|L:\ci L|=N^{n/2}$ and so
$L[\ci]=N^{-1/2}\ci L$ has volume~1. Thus $L[\ci]$ is
a unimodular lattice.

We finally show that $L[\ci]$ is an even unimodular lattice.
Since $L[\ci]$ is integral, to show that it is even it suffices
to show that each vector $\bu$ in a generating set of $L[\ci]$ has
$|\bu|^2$ even. The vectors $\bu=N^{-1/2}\al\bv$ for $\al\in\ci$
and $\bv\in L$ generate $L[\ci]$. Then
$$|\bu|^2=\frac{1}{N}|\al\bv|^2=\frac{|\al|^2}{N}|\bv|^2.$$
But $|\al^2|=\al\al^*\in\ci\ci^*=N\co$ and so $|\al|^2/N\in\Q\cap\co=\Z$
and $|\bv|^2$ is an even integer, as $\bv\in L$, an even lattice.
Thus $|\bu|^2$ is an even integer. Thus $L[\ci]$ is an even
unimodular lattice.

Now suppose that $\ci$ and $\cj$ lie in the same ideal class of~$\co$.
Then $\cj=\al\ci$ where $\al$ is a nonzero element of~$K$. Then
$\cj L=\al\cj_L$ and so
$$L[\cj]=\frac{1}{\sqrt{N(\cj)}}\cj L=\frac{1}{\sqrt{N(\cj)}}\al\ci L
=\sqrt{\frac{N(\ci)}{N(\cj)}}\al L[\ci].$$
Let $\ga=\al\sqrt{N(\ci)/N(\cj)}$. Since $\cj=\al\ci$ then
$N(\cj)=|\al|^2N(\ci)$ and so $|\ga|=1$. By Lemma~\ref{length},
the map $\bv\mapsto\ga\bv$ is an isometry of $\R^n$ and as $L[\cj]=\ga L[\ci]$,
the lattices $L[\ci]$ and $L[\cj]$ are isometric.
\end{pf}

Given~$L$, we can produce a maximum of $h$ non-isometric lattices $L[\ci]$
where $h$ denotes the class-number of the quadratic field~$K$.

It is useful to note which for which ideals $\ci$ is $\ci L_+\sub\Z^n$.

\begin{lem}
Let $\ci$ be an ideal of~$\co$. Then $\ci L_+\sub\Z^n$ if and only if
$\ci\sub\id{2,\frac12(1-\sqrt{-l})}$. In this case also
$N(\ci)\Z^n\sub\ci L_+$.
\end{lem}
\begin{pf}
Note that $L_+\cap\Z^n=L_0$ and so $\ci L_+\sub\Z^n$
if and only if $\ci L_+\sub L_0$. This occurs if and only if
$\ci$ annihilates the $\co$-module $M=L_+/L_0$. This module has
2 elements, so it must be isomorphic to $\co/\cj$ where
$\cj$ is an ideal of norm~2. As $\id{2,\frac12(1-\sqrt{-l})}$
has norm 2 and is seen to annihilate $M$ as $\frac12(1-\sqrt{-l})$
takes $\frac12\be_1(I+W)$ to $\frac14(l+1)\be_1$, then
$\cj=\id{2,\frac12(1-\sqrt{-l})}$. Thus $\cj$ is the annihilator
of $M$ and the first statement follows.

Suppose that $\ci\sub\id{2,\frac12(1-\sqrt{-l})}$. The lattice
$L_+[\ci]$ is unimodular so that if $\bu\cdot\bv\in\Z$ for all
$\bv\in L_+[\ci]$ then $\bu\in L_+[\ci]$. If $\bu=\sqrt{N(\ci)}\bw$
with $\bw\in\Z^n$ then $\bu\cdot\bv\in\Z$ for all $\bv\in N(\ci)^{-1/2}\Z^n$
and as $L_+[\ci]\sub N(\ci)^{-1/2}\Z^n$ then $\bu\in L_+[\ci]$.
Hence $\sqrt{N(\ci)}\Z^n\sub L_+[\ci]$ and so $N(\ci)\Z^n\sub\ci L_+$.
\end{pf}

In this case the lattice $\La$ is the inverse image of a subgroup
$\cc$ of $(\Z/N\Z)^n$, where $N=N(\ci)$,
under the projection $\pi:\Z^n\to(\Z/N\Z)^n$.
Such a subgroup is called a \emph{linear code} of length $n$
over $\Z/N\Z$. We also say that $\La$ is obtained from $\cc$
by \emph{construction A${}_N$}.

The standard dot product is well-defined on the group $(\Z/N\Z)^n$.
If a subgroup $\cc\sub(\Z/N\Z)^n$ satisfies $\bu\cdot\bv=0$
for all $\bu$, $\bv\in\cc$ then $\cc$ is \emph{self-orthogonal}.
Also $\cc$ is \emph{self-dual} if $\bu\cdot\cc=0$ if and only
if $\bu\in\cc$. By the nondegeneracy of the dot product,
$\cc$ is self-dual if and only if $\cc$ is self-orthogonal
and $|\cc|=N^{n/2}$.


\begin{prop}\label{codegens}
Let $\ci$ be an ideal of $\co$ with
$\ci\sub\id{2,\frac12(1-\sqrt{-l})}$ and $N(\ci)=N$. The lattice $L=\ci L_+$
is obtained from construction A${}_N$ from a self-dual linear code
$\cc$ of length $n$ over $\Z/N\Z$.

If $\ci=\id{N,\frac12(a-\sqrt{-l})}$ with $a\equiv1$ (mod~$4$) and
$a^2\equiv-l$ (mod~$4N$) then $\cc$ is spanned by the vectors of the form
$\frac12(\be_i+\be_j)(aI-W)$ ($1\le i\le j\le n$).
\end{prop}
\begin{pf}
Apart from the self-duality of $\cc$ we have already proved the first
assertion. The self-duality of $\cc$ follows from the unimodularity of
$N^{-1/2}\ci L_+$. By volume considerations
$$N^{n/2}=|\Z^n:\ci L_+|=|(\Z/N\Z)^n:\cc|$$
and so $|\cc|=N^{n/2}$. Also if $\bu$, $\bv\in\ci L_+$
then $N^{-1/2}\bu$ and $N^{-1/2}\bv$ lie in the integral
lattice $N^{-1/2}\ci L_+$ so that $N^{-1}\bu\cdot\bv\in\Z$.
Hence $\cc$ is self-orthogonal, and as it has the correct order,
it is self-dual.

The ideal $\ci$ contains the subgroup $N\Z+\frac12(a-\sqrt{-l})\Z$ of $\co$
and as this subgroup also has index $N$ in $\co$ then
$\ci=N\Z+\frac12(a-\sqrt{-l})\Z$. It follows that
$\ci L_+=NL_++\frac12(a-\sqrt{-l})L_+$. As $a\equiv1$ (mod~4),
$\frac12(a+\sqrt{-l})-\frac12(1+\sqrt{-l})$ is an even integer.
It follows that $L_0+\frac12\be_1(aI+W)=L_0+\frac12\be_1(I+W)$
and so $L_+$ is generated by $L_0$ and $\bu=\frac12\be_1(aI+W)$.
Thus $NL_+$ is generated by the $N(\be_i+\be_j)$ and $N\bu$
and $\frac12(a-\sqrt{-l})L_+$ is generated by the $\frac12(\be_i+\be_j)(aI-W)$
and
$$\frac12\bu(aI-W)=\frac14\be_1(aI-W)(aI+W)=\frac{a^2+l}{4}.$$
Note that $(a^2+l)/4$ is a multiple of~$N$. It follows that
$\cc$ is generated by $N\bu$ and the $\frac12(\be_i+\be_j)(aI-W)$.
But $N\bu=(N/2)\be_1(I+W)$ is congruent modulo $N$ to the word
consisting of all $N/2$s. Also $(N/2)\be_1(aI-W)$ is congruent to the
same word. We can drop the generator $N\bu$ and deduce that
$\cc$ is generated by the $\frac12(\be_i+\be_j)(aI-W)$.
\end{pf}

\section{Quadratic residue codes}

To use the above construction of lattices, we need a supply of
skew-symmetric conference matrices. Paley \cite{Pal} constructed a family
of such matrices of order $n=l+1$ whenever $l\equiv3$ (mod~4) is
prime. To apply our theory we stipulate in addition that
$l\equiv7$ (mod~8). We find that the lattices $\ci L_+$ are
derived from quadratic residue codes in this case.

We define a conference matrix $W\in\cw_n$ as follows. Let
$$W=\left(\begin{array}{rrrr}
0&1&\cdots&1\\
-1&&&\\
\vdots&&W'&\\
-1&&&
\end{array}\right)$$
where the $l$-by-$l$ matrix $W'$ is the circulant matrix whose $(i,j)$-entry
is
$$W_{ij}'=\leg{j-i}{l}$$
where $\left(\frac{\ \ }{\ \ }\right)$ denotes the Legendre symbol.
This matrix $W$ is called a conference matrix of \emph{Paley type}.
For the rest of this section $W$ will denote this particular matrix.

We follow the usual practice with quadratic residue codes and label
the entries of a typical vector of length $n=l+1$ using the elements of
the projective line over $\F_l$ as follows:
$\bv=(v_\infty,v_0,v_1,v_2,\ldots,v_{l-1})$. We let
$\be_\infty$, $\be_0$, $\be_1,\ldots\be_{l-1}$
denote the corresponding unit vectors, that is, $\be_\mu$ has a one in the
position labelled $\mu$, and zeros elsewhere.

\begin{lem}[Paley]\label{Paley}
The matrix $W$ is a skew-symmetric conference matrix.
\end{lem}
\begin{pf}
See for instance \cite[Chapter 18]{vLW}.
%
%
\end{pf}

In $\co$, the ideal $2\co$ splits as a product of two distinct prime
ideals: $2\co=\cp\cq$ where $\cp=\id{2,\frac12(1+\sqrt{-l})}$
and $\cq=\cp^*=\id{2,\frac12(1-\sqrt{-l})}$. We shall investigate
the lattices $L_+[\cp^r]$ and $L_+[\cq^r]$ for integers $r\ge0$.
(The discussion for $L_-[\cp^r]$ and $L_-[\cq^r]$ is similar.)

We first need a lemma on the structure of the ideals $\cp^r$
and $\cq^r$.
\begin{lem}\label{powers}
Let $r$ be a positive integer. Then
$$\cp^r=2^r\Z+\frac{1}{2}(t+\sqrt{-l})\Z$$
and
$$\cq^r=2^r\Z+\frac{1}{2}(t-\sqrt{-l})\Z$$
where $t$ is any integer with $t^2\equiv-l$ (mod~$2^{r+2})$ and
$t\equiv1$ (mod~4).
\end{lem}
\begin{pf}
It is well-known that if $s\ge3$, and $a\equiv1$ (mod~8) then the
congruence $x^2\equiv a$ (mod~$2^s$) is soluble. Thus there exists
$t$ with $t^2\equiv-l$ (mod~$2^{r+2}$). By replacing $t$ by $-t$
if necessary, we may assume that $t\equiv1$ (mod~4). Consider the ideal
$\ci=\id{2^r,\frac12(t+\sqrt{-l})}$ of~$\co$. As $2^r\in\ci$
then $\ci$ is a factor of $2^r\co=\cp^r\cq^r$. But as
$\frac12(t+\sqrt{-l})=\frac12(1+\sqrt{-l})+2(t-1)/4\in\cp$
then $\cp$ is a factor of~$\ci$. But $\frac14(t+\sqrt{-l})\notin\co$,
and so $2\co_K=\cp\cq$ is not a factor of~$\ci$. Hence $\ci=\cp^{r'}$
where $1\le r'\le r$. Letting $\al=\frac12(t+\sqrt{-l})$ we have
\begin{eqnarray*}
\ci\ci^*&=&\id{2^r,\al}\id{2^r,\al^*}\\
&=&\id{2^{2r},2^r\al,2^r\al^*,\al\al^*}\\
&=&\id{2^{2r},2^r\al,2^r\al^*,(t^2+l)/4}\\
&\sub&2^r\co
\end{eqnarray*}
as $t^2\equiv-l$ (mod~$2^{r+2}$). But $\ci\ci^*=N(\ci)\co=2^{r'}\co$
and so $r=r'$, that is $\ci=\cp^r$.

Now $\cp^r\sub2^r\Z+\frac12(t+\sqrt{-l})\Z$, but both these groups have index
$2^r$ in~$\co$ so they are equal. The statement about $\cq^r$ now
follows by complex conjugation.
\end{pf}

%

We now consider the lattices $\cq^r L_+$ for $r\ge1$.
Since $\cq^r\sub\cq$ and $\cq=\id{2,\frac12(1-\sqrt{-l})}$
then by Proposition~\ref{codegens}
$\cq^r L_+$ is obtained by construction A${}_{2^r}$ from a self-dual code
$\cc_r$ over $\Z/2^r\Z$.
We shall show that $\cc_r$ is the Hensel lift of an extended
quadratic residue code in the sense of~\cite{BSC}.

%

Recall that the integer $t$ satisfies $t\equiv1$ (mod~4) and $t^2\equiv-l$
(mod~$2^{r+2}$). By Proposition~\ref{codegens} it follows that $\cc_r$
is generated by the vectors $\frac12(\be_i+\be_j)(W-tI)$. It is plain that
we need only these vectors with $i=\infty$ and so $\cc_r$ is spanned by
$\bu=\be_\infty(W-tI)$ and $\bv_j=\frac12(\be_\infty+\be_j)(W-tI)$
for $0\le j<l$.

%
%
%

Let $\phi:(\Z/2^r\Z)^n\to(\Z/2^r\Z)^l$ be the map given by deleting the
first coordinate of the vector. The code $\cc_r$ contains the vector
$\bu=(-t,1,1,\ldots,1)$. As $r$ is odd and $\cc_r$ is self-dual,
the intersection of $\cc_r$ and the kernel of $\phi$ is trivial.
Thus $\cc_r'=\phi(\cc_r)$ has the same order as~$\cc_r$.
Then $\phi(\bu)$ is the all-ones vector, and $\phi(\bv_j)$ are cyclic
shifts of $\phi(\bv_0)$. Also $\phi(\bv_0)=(c_0,c_1,\ldots,c_{p-1})$
where
$$c_j=\left\{\begin{array}{cl}
(1-t)/2&\textrm{if $j=0$,}\\
1&\textrm{if $j$ is a quadratic residue modulo~$l$,}\\
0&\textrm{if $j$ is a quadratic nonresidue modulo~$l$.}
\end{array}\right.$$
Thus $\cc_r'$ is a cyclic code over $\Z/2^r\Z$.

We recall the definition of quadratic residue codes. Consider the
polynomial $X^l-1$ over the field $\F_2=\Z/2\Z$. Then $X^l-1$ splits
into linear factors in some finite extension $\F_{2^k}$ of~$\F_2$.
In fact
$$X^l-1=\prod_{j=0}^{l-1}(X-\ze^j)$$
where $\ze$ is a primitive $l$-th root of unity in~$\F_{2^k}$.
We write
$$X^l-1=(X-1)f_+(X)f_-(X)$$
where
$$f_+(X)=\prod_{(j/l)=1}(X-\ze^j)\quad\textrm{and}
\quad f_-(X)=\prod_{(j/l)=-1}(X-\ze^j).$$
As $l\equiv7$ (mod~8), then 2 is a quadratic residue modulo~$l$, and
so the coefficients of both $f_+$ and $f_-$ are invariant under the
Frobenius automorphism $\de\mapsto\de^2$ of~$\F_{2^k}$. Consequently
both $f_+$ and $f_-$ have coefficients in~$\F_2$. The labelling of
these factors as $f_+$ and $f_-$ depends on the choice of~$\ze$.
Replacing $\ze$ by another primitive $l$-th root of unity either
preserves or interchanges $f_+$ and $f_-$. The coefficients of
$X^{(l-3)/2}$ in $f_+$ and $f_-$ are 0 and 1 is some order, so we can,
and shall, choose $\ze$ such that
$$f_+(X)=X^{(l-1)/2}+0X^{(l-3)/2}+\cdots\quad\textrm{and}
\quad f_-(X)=X^{(l-1)/2}+X^{(l-3)/2}+\cdots.$$
The cyclic codes of length $l$ over $\F_2$ with generator polynomials
$f_+(X)$ and $f_-(X)$ are called the quadratic residue codes.

Bonnecaze, Sol\'e and Calderbank \cite{BSC}
extended the notion of quadratic residue code to codes over
$\Z/2^r\Z$. By Hensel's lemma there exist unique polynomials
$f_+^{(r)}(X)$ and $f_-^{(r)}(X)$ with coefficients in $\Z/2^r\Z$ such that
$$X^l-1=(X-1)f_+^{(r)}(X)f_-^{(r)}(X),$$
$$f_+^{(r)}(X)\equiv f_+(X)\pmod{2}\quad\textrm{and}\quad
f_-^{(r)}(X)\equiv f_-(X)\pmod{2}.$$
The cyclic codes over $\Z^l$ with generator polynomials
$f_+^{(r)}(X)$ and $f_-^{(r)}(X)$ are called
lifted quadratic residue codes over $\Z/2^r\Z$.

\begin{thm}\label{QRcode}
The code $\cc_r'$ is the lifted quadratic residue code over $\Z/2^r\Z$
with generator polynomial~$f_+^{(r)}(X)$.
\end{thm}
\begin{pf}
Cyclic codes of length $l$ over $R=\Z/2^r\Z$ correspond to ideals of the
polynomial ring $R[X]/\id{X^l-1}$. The code $\cc_r'$ corresponds to the
ideal $\ci=\id{g,h}$ where
$$g(X)=\sum_{j=0}^{p-1}X^j$$
and
$$h(X)=\frac{1-t}{2}+\sum_{(j/l)=1}X^j.$$

We first consider the case where $r=1$. Then $\ci=\id{u(X)}$ where $u(X)$
is the greatest common divisor of $g(X)$ and~$h(X)$. Let $\ze$ be a root
of $f_+(X)=0$ in an extension field of~$\F_2$. The roots of $g(X)$
are the $\ze^j$ where $p\nmid j$. As $t\equiv1$ (mod~4) then $\frac12(1-t)$
is even and so $h(X)=\sum_{(j/l)=1}X^j$. Now
$$\sum_{(j/l)=1}\ze^j=0\qquad\textrm{and}\qquad\sum_{(j/l)=-1}\ze^j=1.$$
It follows that
$$h(\ze^a)=\sum_{(j/l)=1}(\ze^a)^j=0$$
if and only if $\leg{a}{l}=1$. Thus $u(X)=f_+(X)$.

Now we consider the general case. The reduction of $\cc_r'$ modulo 2
is~$\cc_1'$. Any liftings to $\cc_r'$ of a basis of~$\cc_1'$ generate
a free $R$-module (of rank~$\frac12(l+1)$), and so they generate the whole
code~$\cc_r'$. As $\cc_r$ is free over~$R$, it is generated as an ideal
by a monic polynomial~$F(X)$, of degree $\frac12(l+1)$. As $F(X)$ reduces
to $f_+(X)$ modulo~2, and $F(X)\mid X^l-1$ it follows that $F(X)=f_+^{(r)}(X)$
as required.
\end{pf}

Given the code $\cc_r'$, the code $\cc_r$ can be reconstructed, since
for each element of $\cc_r'$ the corresponding element of $\cc_r$ is uniquely
determined as it is orthogonal to $(-t,1,1,\ldots,1)$.

We now turn to $\cp^r L_+$. This is no longer a sublattice of~$\Z^n$.

\begin{lem}\label{generators2}
Let $r$ be a positive integer. The index $|\cp^r L_+:\cp^r L_+\cap\Z^n|=2$.
The lattice $\cp^r L_+\cap\Z^n$ is generated by the vectors
$2^r(\be_\infty+\be_\mu)$ ($\mu\in\{\infty,0,1,2,\ldots,l-1\}$),
the vector $\bu=\be_\infty(W+tI)$ and
and the vectors $\bv_j=\frac12(\be_\infty+\be_j)(W+tI)$ ($0\le j<l$).
Also $\cp^r L_+$ is generated by
$\cp^r L_+\cap\Z^n$ and $\frac12\bu-\frac14(t^2+l)\be_\infty$.
\end{lem}
\begin{pf}
We have $\cp^r=2^r\Z+\frac12(t+\sqrt{-l})\Z$. Let $\Om_0$ denote the
lattice generated by the $2^r(\be_\infty+\be_\mu)$, $\bu$ and the~$\bv_j$.
The lattice $L_+$ is generated
by the $\be_\infty+\be_\mu$ and $\frac12\be_\infty(tI+W)$. Thus
$2^r(\be_\infty+\be_\mu)$, $\bu=\frac12(t+\sqrt{-l})2\be_\infty$ and
$\bv_j=\frac12(t+\sqrt{-l})(\be_\infty+\be_j)$ all lie in $\cp^rL_+$.
These vectors all have integer coordinates, and so
$\Om_0\sub\cp^r L_+\cap\Z^n$.

Let $\Om$ be the lattice generated by $\Om_0$ and
$\frac12\bu-\frac14(t^2+l)\be_\infty$. Now
\begin{eqnarray*}
\frac12(t+\sqrt{-l})\frac12\be_\infty(tI+W)
&=&\frac14(t+\sqrt{-l})^2\be_\infty\\
&=&\left[\frac{t}{2}(t+\sqrt{-l})
-\frac{t^2+l}{4}\right]\be_\infty\\
&=&\frac{t}{2}\bu-\frac{t^2+l}{4}\be_\infty.
\end{eqnarray*}
As $t$ is odd and $\bu\in\Om_0$ then $\Om\sub\cp^r L_+$.

The lattice $\cp^r L_+$ is generated by $\Om$ and $2^{r-1}\be_\infty(tI+W)
=2^{r-1}\bu$. But $\bu-\frac12(t^2+l)\be_\infty\in\Om$ and as $t^2+l$
is divisible by $2^{r+1}$ then $\bu\in\Om_0$ and so $\Om=\cp^r L_+$.
As $\frac12\bu-\frac14(t^2+l)\be_\infty$ is not in $\Z^n$ but its double
is in $\Om_0$, then $|\Om:\Om_0|=|\cp^r L_+:\cp^r L_+\cap\Z^n|=2$
and so $\Om_0=\cp^r L_+\cap\Z^n$.
\end{pf}

One can now proceed to express the lattices $\cp^r L_+$ and $\cp^r L_+\cap\Z^n$
in terms of lifted quadratic residue codes over $\Z/2^r\Z$. For simplicity
we present the details only for $r=1$. Let $\cd'$ denote the cyclic quadratic
residue code of length $l$ over $\F_2$ with generator polynomial $f_-(X)$,
and let $\cd$ denote its extension obtained by appending a parity check bit
at the front.

\begin{thm}\label{QRcode2}
The lattice $\cp L_+\cap\Z^n$ consists of those vectors reducing modulo
2 to elements of $\cd$ and the sum of whose entries is a multiple of 4.
The lattice $\cp L_+$ is obtained from $\cp L_+\cap\Z^n$ by adjoining the
extra generator $\frac12(\frac12(1-l),1,1,\ldots,1)$.
\end{thm}
\begin{pf}
We may take $t=1$ in the proof of Lemma~\ref{generators2}. In this case
the vector $\bu$ is the all-ones vector while each $\bv_j$ consists of
$\frac12(l+1)$ ones and $\frac12(l+1)$ zeros. As $\frac12(l+1)$
is a multiple of 4 the sum of the entries of each of these vectors
is a multiple of 4. As this is manifestly true for the vectors
$2(\be_\infty+\be_\mu)$ too, then the sum of the entries of each vector in
$\cp L_+\cap\Z^n$ is a multiple of~4.

If we delete the first entry of the given generators of $\cp L_+$
and reduce modulo 2 we get the all-ones vector of length $l$ and the
cyclic shifts of the vector $\bw_0=(d_0,d_1,\ldots,d_{l-1})$
where
$$d_j=\left\{\begin{array}{cl}
1&\textrm{if $j=0$ or if $j$ is a quadratic residue modulo~$l$,}\\
0&\textrm{if $j$ is a quadratic nonresidue modulo~$l$.}
\end{array}\right.$$
By a similar argument to the proof of Theorem~\ref{QRcode} these
vectors generate the cyclic quadratic residue code~$\cd'$.
Hence each element of $\cp L_+\cap\Z^n$ reduces modulo 2 to an element
of~$\cd$. If $\Om$ denotes the sublattice of $\Z^n$ consisting of
vectors reducing modulo 2 to $\cd$ and with the entries summing
to a multiple of~4, then $|\Z^n:\Om|=2^{1+n/2}=|\Z^n:\cp L_+\cap\Z^n|$.
Thus $\Om=\cp L_+\cap\Z^n$.

Now letting $t=1$ we see that $\frac12\bu-\frac14(t^2+l)\be_\infty
=\frac12(\frac12(1-l),1,1,\ldots,1)$ and so this vector together
with $\cp L_+\cap\Z^n$ generates~$\cp L_+$.
\end{pf}

In the terminology of Conway and Sloane \cite[Chapter~5, \S3]{CS},
the lattice $\cp L_+\cap\Z^n$
is obtained from the code $\cd$ by construction B. Then the lattice
$\cp L_+$ is obtained by density doubling. One can extend these
notions to lifted quadratic residue codes to produce the lattices
$\cp^r L_+$.

We look briefly at the lattices $\ci L_+$ for more general ideals~$\ci$.
Consider the case where $\ci=\ca$, an ideal of norm $p$, an odd prime.
Then $\ca=\id{p,t+\sqrt{-l}}$ where $t^2\equiv-l$ (mod~$p$).
The rows of the matrix $tI+W$ generate a self-dual linear code $\cc$ over
$\F_p$ which turns out to be an extended quadratic residue code.
The coordinates of vectors in $L_+$ are half-integers, and it is meaningful
to reduce these modulo the odd prime~$p$. Then the lattice $\ca L_+$
simply consists of the vectors in $L_+$ which reduce modulo $p$ to
elements of~$\cc$. More generally $\ca^r L_+$ will have a similar description
in terms of an extended lifted quadratic residue code over $\Z/p^r\Z$. Finally
by splitting a general ideal $\ci$ into a product of powers of prime
ideals~$\ca^r$, we can describe $\ci L_+$ in terms of the various
$\ca^r$ using the Chinese remainder theorem.

\section{Examples}

Since the ring $\Z[\frac12(1+\sqrt{-7})]$ has class number 1 (and each
even unimodular rank 8 lattice is isometric to the $E_8$ root lattice)
the first interesting examples occur when $l=15$ and the
first interesting examples with Paley matrices occur when $l=23$.

\subsection{$l=23$ and $l=31$}

In both these cases we take $W$ to be the Paley matrix. We first consider
the case $l=23$.

The class group of $\co=\Z[\frac12(1+\sqrt{-23})]$ has order~3, and
the class of each of its ideals $\cp=\id{2,\frac12(1+\sqrt{-23})}$ and
$\cq=\id{2,\frac12(1-\sqrt{-23})}$ generates its class group.
The lattice $L_+$ itself is the lattice $D_{24}^+$. The lattices
$\cq^r L_+$ are obtained by applying construction A${}_{2^r}$ to the lifted
quadratic residue codes~$\cc_r$. The code $\cc_1$ is the extended
binary Golay code. It is plain that $\cq L_+$ is obtained by applying
construction A [Chapter 5, \S2] to the binary Golay code, and so $L_+[\cq]$
is isometric to the Niemeier lattice with root system~$A_1^{24}$.

The isometry classes of the unimodular
lattices $L_+[\cq^r]$ depend only on the congruence class of $r$
modulo~3. If $r\equiv0$ (mod~3) then $L_+[\cq^r]$ is isometric to
$D_{24}^+$ while if $r\equiv1$ (mod~3) then $L_+[\cq^r]$ is isometric to
the Niemeier lattice with root system~$A_1^{24}$. To identify $L_+[\cq^r]$
when $r\equiv2$ (mod~3) note that $\cq^2$ lies in the same ideal class
as~$\cp$. Hence for $r\equiv2$ (mod~3), $L_+[\cq^r]$ is isometric
to~$L_+[\cp]$. By Theorem~\ref{QRcode2} it is plain that
$L_+[\cp]$ is the Leech lattice, as given by Leech's original
construction~\cite{Leech}. Applying Theorem~\ref{main} gives
an explicit isomorphism between $L_+[\cp]$ and $L_+[\cq^2]$ which
is equivalent to that constructed in \cite{ChS}.

In general if $s$ is the order of the class of the ideal $\cp$
in the class group of~$\co$, then up to isometry $L_+[\cp^r]$ and
$L_+[\cq^r]$ depend only on the congruence class of $r$ modulo~$s$. Also
$L_+[\cp^r]$ and $L_+[\cq^{r'}]$ will be isometric whenever
$r\equiv -r'$ (mod~$s$). For $l=31$ we also have $s=3$ and the above
discussion is valid for $l=31$ too. In particular
$L_+[\cp]$ is isometric to $L_+[\cq^2]$,
and we recover \cite[Theorem~1]{ChS}.

We can give alternative constructions of the Leech lattice at will
simply by writing down ideals of $\Z[\frac12(1+\sqrt{-23})]$
equivalent to~$\cp$. Let $\ci=\id{3,\frac12(1+\sqrt{-23})}$
and $\cj=\id{3,\frac12(-1+\sqrt{-23})}$. Then $\cp$,
$\cj$ and $\cq\ci=\id{6,\frac12(-5+\sqrt{-23})}$ all
lie in the same ideal class.

The lattice $\cj L_+$ is generated using density doubling from the
lattice $L'$ consisting of all vectors in $\Z^{24}$ whose
entries sum to zero and which reduce modulo 3 to elements of the
extended ternary quadratic residue code with generator matrix $I-W$.
Then $\cj L_+$
is generated by $L'$ and the vector $\frac12(5,1,1,\ldots,1)$.
The lattice $L_+[\cj]=3^{-1/2}\cj L_+$ is isometric to the
Leech lattice.

Next consider the lattice $\cq\ci L_+$. This consists of the vectors
in $\Z^{24}$ reducing modulo 2 and modulo 3 to elements of appropriately
chosen binary and ternary quadratic residue codes. The binary code
in question is that generated by vectors $\frac12(\be_\infty+\be_\al)(I-W)$
for $\al\in\{\infty,0,1,2,\ldots,l-1\}$ and the ternary code
is generated by the rows of~$I+W$. Then  $L_+[\cq\ci]=6^{-1/2}\cq\ci L_+$
is isometric to the Leech lattice.

\subsection{$l=47$}

Again we take $W$ to be the Paley matrix.
In \cite[Chapter 7, \S7]{CS} the lattice $\La=P_{48q}$ is described.
This is an even unimodular lattice of rank 48 and minimum norm~6.
It is generated
by the following vectors $(a_\infty,a_0,a_1,\ldots,a_{46})$:
\begin{description}
\item{(i)} $(1/\sqrt{12})(-5,1,1,\ldots,1)$,
\item{(ii)} those vectors of the shape $(1/\sqrt{3})(1^{24},0^{24})$
supported on the translates modulo 47 of the set $\{0\}\cup Q$
where $Q$ is the set of quadratic residues modulo 47,
\item{(iii)} all vectors of the shape $(1/\sqrt{3})(\pm3^2,0^{46})$.
\end{description}
It is more convenient to consider instead the equivalent lattice $\La'$
generated by the vectors
\begin{description}
\item{(i)${}'$} $(1/\sqrt{12})(5,1,1,\ldots,1)$,
\item{(ii)${}'$} those vectors of the shape $(1/\sqrt{3})(1^{24},0^{24})$
supported on the translates modulo 47 of the set $\{0\}\cup N$
where $N$ is the set of quadratic nonresidues modulo 47,
\item{(iii)${}'$} all vectors of the shape $(1/\sqrt{3})(\pm3^2,0^{46})$.
\end{description}

We claim that $\La'$ is the lattice $L_+[\ci]$ where
$\ci=\id{3,\frac12(1-\sqrt{-47})}$. Note that the norm of $\ci$
is~3. It suffices to show that each of
the generating vectors for $\La'$ is contained in~$L_+[\ci]$. Since
each vector of shape $(\pm1^2,0^{46})$ lies in $L_+$ and $3\in P$
then it is immediate that the vectors of type (iii)${}'$ lie in $L_+[\ci]$.
The vectors of type (ii)${}'$ are the differences of the first row
and an arbitrary other row of the matrix $(1/2\sqrt{3})(I-W)$.
Since $\frac12(1-\sqrt{-47})\in\ci$, the vectors of type (ii)${}'$
lie in~$L_+[\ci]$. Finally, $\frac12(1,-1,-1,\ldots,-1)$, the
first row of $\frac12(I-W)$, lies in $\ci L_+$. Also $\bv_0=\frac12\be_0(I+W)
\in L_+$ and $3\bv_0=\frac12(3,3,3,\ldots,3)\in\ci L_+$. Adding these
two vectors gives $\frac12(5,1,1,\ldots,1)\in\ci L_+$ so that the vector
of type (i)${}'$ does lie in~$L_+[\ci]$.

The ideal $\id{\frac12(1-\sqrt{-47})}$ has norm 12 and factors
as $\cq^2\ci$. The class number of $\Q(\sqrt{-47})$ is 5, and so $[\ci]=[\cp^2]
=[\cq^3]$. Thus $\La'$ is isometric to $L_+[\cq^3]$, which is constructed
using construction A from the quadratic residue code of length 48
over $\Z/8\Z$.

\subsection{$l=15$}

In this case there is no Paley matrix. We consider two different
conference matrices of order~16.

If $W\in\cw_n$ then the $2n$-by-$2n$ matrix
$$W'=\left(\begin{array}{cc}
W&I+W\\
-I+W&-W
\end{array}\right)$$
is a skew-symmetric conference matrix of order~$2n$. Applying this
construction four times to the zero matrix in $\cw_1$ gives the
matrix
$$W_1=\left(\begin{array}{cccccccccccccccc}
0&+&+&+&+&+&+&+&+&+&+&+&+&+&+&+\\
-&0&-&+&-&+&-&+&-&+&-&+&-&+&-&+\\
-&+&0&-&-&+&+&-&-&+&+&-&-&+&+&-\\
-&-&+&0&-&-&+&+&-&-&+&+&-&-&+&+\\
-&+&+&+&0&-&-&-&-&+&+&+&+&-&-&-\\
-&-&-&+&+&0&+&-&-&-&-&+&+&+&+&-\\
-&+&-&-&+&-&0&+&-&+&-&-&+&-&+&+\\
-&-&+&-&+&+&-&0&-&-&+&-&+&+&-&+\\
-&+&+&+&+&+&+&+&0&-&-&-&-&-&-&-\\
-&-&-&+&-&+&-&+&+&0&+&-&+&-&+&-\\
-&+&-&-&-&+&+&-&+&-&0&+&+&-&-&+\\
-&-&+&-&-&-&+&+&+&+&-&0&+&+&-&-\\
-&+&+&+&-&-&-&-&+&-&-&-&0&+&+&+\\
-&-&-&+&+&-&+&-&+&+&+&-&-&0&-&+\\
-&+&-&-&+&-&-&+&+&-&+&+&-&+&0&-\\
-&-&+&-&+&+&-&-&+&+&-&+&-&-&+&0
\end{array}\right)$$
where, for convenience, we have denoted $1$ and $-1$ by $+$ and $-$
respectively. The ideal class group of $\Z[\frac12(1+\sqrt{-15})]$
has order~2. The ideal $\ci=\id{2,\frac12(1-\sqrt{-15})}$
is not principal and $\ci L^+$ is given by construction A from the
binary code with the generator matrix
$$\left(\begin{array}{cccccccccccccccc}
1&0&0&1&0&1&0&1&0&0&0&0&0&0&0&0\\
0&1&0&1&0&0&1&1&0&0&0&0&0&0&0&0\\
0&0&1&1&0&1&1&0&0&0&0&0&0&0&0&0\\
0&0&0&0&1&1&1&1&0&0&0&0&0&0&0&0\\
0&0&0&0&0&0&0&0&1&0&0&0&0&1&1&1\\
0&0&0&0&0&0&0&0&0&1&0&0&1&0&1&1\\
0&0&0&0&0&0&0&0&0&0&1&0&1&1&0&1\\
0&0&0&0&0&0&0&0&0&0&0&1&1&1&1&0
\end{array}\right).$$
Thus $L_+[\ci]$ is isometric to the orthogonal direct sum of two
copies of the $D_8^+$ lattice. This is not isometric to $L_+$.

Another conference matrix of order 16 is
$$W_2=\left(\begin{array}{cccccccccccccccc}
0&+&+&+&+&+&+&+&+&+&+&+&+&+&+&+\\
-&0&+&+&-&+&-&-&+&-&+&+&-&+&-&-\\
-&-&0&+&+&-&+&-&+&-&-&+&+&-&+&-\\
-&-&-&0&+&+&-&+&+&-&-&-&+&+&-&+\\
-&+&-&-&0&+&+&-&+&+&-&-&-&+&+&-\\
-&-&+&-&-&0&+&+&+&-&+&-&-&-&+&+\\
-&+&-&+&-&-&0&+&+&+&-&+&-&-&-&+\\
-&+&+&-&+&-&-&0&+&+&+&-&+&-&-&-\\
-&-&-&-&-&-&-&-&0&+&+&+&+&+&+&+\\
-&+&+&+&-&+&-&-&-&0&-&-&+&-&+&+\\
-&-&+&+&+&-&+&-&-&+&0&-&-&+&-&+\\
-&-&-&+&+&+&-&+&-&+&+&0&-&-&+&-\\
-&+&-&-&+&+&+&-&-&-&+&+&0&-&-&+\\
-&-&+&-&-&+&+&+&-&+&-&+&+&0&-&-\\
-&+&-&+&-&-&+&+&-&-&+&-&+&+&0&-\\
-&+&+&-&+&-&-&+&-&-&-&+&-&+&+&0
\end{array}\right).$$
In this case the lattice $\ci L_+$ is obtained using construction A
applied to the binary code with generator matrix
$$\left(\begin{array}{cccccccccccccccc}
1&0&0&0&0&0&0&0&0&1&1&1&1&1&1&1\\
0&1&0&0&0&0&0&0&1&0&1&1&1&1&1&1\\
0&0&1&0&0&0&0&0&1&1&0&1&1&1&1&1\\
0&0&0&1&0&0&0&0&1&1&1&0&1&1&1&1\\
0&0&0&0&1&0&0&0&1&1&1&1&0&1&1&1\\
0&0&0&0&0&1&0&0&1&1&1&1&1&0&1&1\\
0&0&0&0&0&0&1&0&1&1&1&1&1&1&0&1\\
0&0&0&0&0&0&0&1&1&1&1&1&1&1&1&0
\end{array}\right).$$
Thus $L_+[\ci]$ is isometric to the $D_{16}^+$ lattice and so to $L_+$.
This example shows that the isometry class of $L_+[\ci]$
depends on the choice of the conference matrix $W$, and also
that $L_+[\ci]$ and $L_+[\cj]$ may be isometric even when $\ci$
and $\cj$ are in different ideal classes.


\end{document}